\ifx\documentclass\undefined
\documentstyle[10pt]{article}
\else
\documentclass[10pt]{article}
\usepackage{latexsym}
\usepackage{amsfonts}
\usepackage{amsmath}
\usepackage{amssymb}
\fi

\author{K. Bezdek\thanks{Partially supported by a Natural Sciences and
Engineering Research Council of Canada Discovery Grant.}
\and A. E. Litvak$^{\dagger}$ }

\date{}

\font\tenBbb=msbm10 at 12pt         \font\sevenBbb=msbm9    \font\fiveBbb=msbm7
\newfam\Bbbfam
\textfont\Bbbfam=\tenBbb \scriptfont\Bbbfam=\sevenBbb
\scriptscriptfont\Bbbfam=\fiveBbb

\def\R{{\mathbb R}}

\def\kkk{\null\hfill $\Box$\smallskip}
\def\r{ \right}

\def\lam{\lambda}

\def\la{\left\langle}
\def\ra{\r\rangle}

\def\KK{{\bf K}}
\def\LL{{\bf L}}

\def\BB{{\bf B}}

\newcommand{\proof}{{\noindent\bf Proof:{\ \ }}}

\newtheorem{theorem}{Theorem}[section]
\newtheorem{lemma}[theorem]{Lemma}
\newtheorem{remark}[theorem]{Remark}

\newtheorem{definition}[theorem]{Definition}

\newcommand{\conv}{\mathop{\rm conv\,}}

\newcommand{\vol}{{\rm vol}}

\newcommand{\crv}{{\rm crv}}

\title{Packing convex bodies by cylinders
}

\begin{document}

\maketitle

\begin{abstract}
In \cite{BL} in relation to the unsolved Bang's plank problem (1951) we obtained a lower bound for the sum of relevant measures of cylinders covering a given $d$-dimensional convex body. In this paper we provide the packing counterpart of these estimates. We also extend bounds to the case of $r$-fold covering and packing and show a packing analog of  Falconer's results
(\cite{Fa}).
\end{abstract}

\section{Introduction}
\label{zero}

In the remarkable paper \cite{Ba} Bang has given an elegant proof of the plank conjecture of
Tarski showing that if a convex body is covered by finitely many planks in $d$-dimensional
Euclidean space, then the sum of the widths of the planks is at least as large as the minimal
width of the body.
We refer to \cite{AKP} for historical remarks and references.
A celebrated extension of Bang's theorem to $d$-dimensional normed spaces
has been given by Ball in \cite{Ball2}. In his paper Bang raises also the important related
question whether the sum of the base areas of finitely many cylinders covering a $3$-dimensional
convex body is at least half of the minimum area of a 2-dimensional projection of the body. In the recent
paper \cite{BL} the authors have investigated this problem of Bang in $d$-dimensional
Euclidean space. In particular, we proved Bang's conjecture with constant one third instead of one half.
From the point of view of discrete geometry it is quite surprising that so far there has not been any
packing analogue of the above theorems on coverings by planks and cylinders. In this paper we fill this gap.

\section{Notation}
\label{one}

We identify a $d$-dimensional affine space with $\R ^d$.
By $|\cdot|$ and $\la \cdot , \cdot \ra$ we denote the canonical
Euclidean norm and the canonical inner product on $\R ^d$.
The canonical Euclidean ball and sphere in $\R^d$ are denoted by
$\BB_2^d$ and $S^{d-1}$. The volume of $\BB_2^d$ is denoted
by $\omega_d$.

By a convex body $\KK$ in $\R^d$ we always mean a compact convex set with  non-empty interior, which is denoted by $\mbox{int} (\KK)$. The volume of a convex body $\KK$ in $\R^d$ is denoted by $\vol (\KK)$. When we would like to emphasize that we take $d$-dimensional volume of a body in $\R^d$ we write $\vol _d (\KK)$.
%
%
%
%
%

The {\it Banach-Mazur distance} between two convex
bodies $\KK$ and $\LL$ in $\R ^d$ is defined by
$$
    d(\KK, \LL) = \inf{\left\{ \lam > 0       \ \mid \ a\in \LL, \ b\in \KK, \
   \LL - a \subset T \left(\KK - b \r) \subset \lam \left(\LL - a\r) \r\}},
$$
where the infimum is taken over all (invertible) linear operators $T : \R^d \to \R^d$.
We denote by $d_{\KK}$ the Banach-Mazur distance between $\KK$ and the Euclidean
ball $\BB_2^d$.
John's Theorem (\cite{J}) implies that for every convex body
$\KK$ in $\R^d$, $d_{\KK}$ is bounded by $d$, moreover if $\KK$
is $0$-symmetric, i.e., symmetric about the origin $0$ in $\R^d$, then $d_{\KK}\leq \sqrt{d}$ (see e.g. \cite{Ball}).



Given a (linear) subspace $E\subset \R^d$ we denote the orthogonal projection
on $E$ by $P_E$ and the orthogonal complement of $E$ by $E^{\perp}$.
We will use the following theorem, proved
by Rogers and Shephard (\cite{RS}, see also \cite{Ch} and Lemma 8.8 in
\cite{Pi2}).

\begin{theorem}\label{rsh}
Let $1\leq k \leq d$. Let $\KK$ be a convex body in $\R^d$ and
$E$ be a $k$-dimensional subspace of $\R^d$. Then
$$
  \max _{x\in \R^d}\ \vol _{d-k} \left( \KK\cap \left( x + E^{\perp}
  \r)\r) \vol _k (P_E\KK) \leq {d \choose k}   \vol _d (\KK).
$$
\end{theorem}

\medskip

\begin{remark}
{\rm Note that the reverse estimate
$$
  \max _{x\in \R^d}\ \vol _{d-k} \left( \KK\cap \left( x + E^{\perp}
  \r)\r) \vol _k (P_E\KK) \geq   \vol _d (\KK)
$$
is a simple consequence of the Fubini Theorem and holds for
every measurable set $\KK$ in $\R^d$.}
\end{remark}


\section{Preliminary results}
\label{two}

Given $0<k<d$ define a $k$-codimensional cylinder $C$ as a set, which can be
presented in the form $C = B + H$, where $H$ is a $k$-dimensional (linear) subspace
of $\R^d$ and  $B$ is a measurable set in $E: =H^{\perp}$.  Given a convex
body $\KK$ and a $k$-co\-di\-men\-si\-o\-nal cylinder $C= B + H$ denote
the crossectional volume of $C$ with respect to $\KK$ by
$$
 \crv _{\KK} (C) := \frac{\vol _{d-k} (C\cap E)}{\vol _{d-k} (P_E
 \KK)} =\frac{\vol _{d-k} (P_E C)}{\vol _{d-k} (P_E
 \KK)} = \frac{\vol _{d-k} (B)}{\vol _{d-k} (P_E \KK)} .
$$
In \cite{BL} (see Remark~2 following Theorem~3.1 there) we proved that if a
convex body $\KK$ is covered by $k$-codimensional cylinders $C_1$, \ldots, $C_N$,
then
\begin{equation}\label{cov-gen}
    \sum _{i=1}^N  \mbox{\rm crv} _{\KK} (C_i) \geq \frac{1}{{d \choose k}} .
\end{equation}
The case $k=d-1$ corresponds to the affine plank problem of Bang (\cite{Ba}), because in this case one has the sum of the relative widths of the planks (i.e., $(d-1)$-codimensional cylinders)
on the left side of (\ref{cov-gen}).
Note that Ball (\cite{Ball2}) proved that such sum should exceed
1 in the case of centrally symmetric convex body $\KK$, while the general case is still open.
The estimate (\ref{cov-gen}) implies the lower bound $1/d$.
Moreover, if $\KK$ is an ellipsoid and $k=1$ one has (see
Theorem~3.1 in \cite{BL})
\begin{equation}\label{cov-ell}
    \sum _{i=1}^N  \mbox{\rm crv} _{\KK} (C_i) \geq 1.
\end{equation}

Akopyan, Karasev and Petrov (\cite{AKP}) have recently proved that (\ref{cov-ell}) holds for $2$-codimensional cylinders
as well. They have also conjectured that (\ref{cov-ell}) holds for $k$-codimensional cylinders for all $0<k<d$.

Before passing to packing, we would like to mention that methods developed in \cite{BL} can be used to prove bounds for multiple coverings. The notion of multiple covering (resp., packing) was introduced in a geometric setting independently by Harold Davenport and L\'aszl\'o Fejes T\'oth \cite{Fe}. Recall that that sets  $\LL _1$, \ldots, $\LL _N$ form an $r$-fold covering of $\KK$ if every point in $\KK$ belongs to at least $r$ of $\LL _i$'s. Slightly modifying proofs of Theorem~1 and Remark~2 in \cite{BL}, we obtain the following theorem.

\begin{theorem}\label{coveringcyl}
Let $\KK$ be a convex body in $\R^d$ and $0<k<d$. Let $C_1, \dots, C_N$
be $k$-co\-di\-men\-si\-o\-nal cylinders in $\R^d$ which form an
$r$-fold covering of $\KK$.
Then
$$
    \sum _{i=1}^N  \mbox{\rm crv} _{\KK} (C_i) \geq  \frac{r}{{d \choose k}} .
$$
Moreover, if $k=1$ and $\KK$ is an ellipsoid then
$$
    \sum _{i=1}^N  \mbox{\rm crv} _{\KK} (C_i) \geq r.
$$
\end{theorem}
\bigskip

\begin{remark}\label{multiple-covering-by-planks}
{\rm The following multiple covering version of Ball's theorem (\cite{Ball2}) seems to be an open problem:
Let $1<r\le N$ be integers. Let $C_1, \dots, C_N$
be planks (i.e., $(d-1)$-co\-di\-men\-si\-o\-nal cylinders) in $\R^d$ which form an
$r$-fold covering of the $0$-symmetric convex body $\KK$. Then prove or disprove that the sum of the relative
widths of the planks is at least $r$, i.e.,
$ \sum _{i=1}^N  \mbox{\rm crv} _{\KK} (C_i) \geq r$. In particular, does the above problem have a positive answer for $r=2$?
We note that Falconer (\cite{Fa}) asked
for a multiple covering version of Bang's theorem and proved such
a result for convex bodies whose minimal width is two times the inradius (including $0$-symmetric convex bodies)
in $\R^2$ and $\R^3$.}
\end{remark}

\bigskip


\section{Packing by cylinders}

In this section we provide estimates for packing by cylinders in terms of
the volumetric parameter, $\crv _{\KK} (C)$,  introduced in \cite{BL}.
Our proofs are close to the proofs of corresponding covering results
in Section~3 of \cite{BL}. We provide all the details for the sake of
completeness. We start with a definition for a packing by cylinders.

\medskip

\begin{definition}
{\rm Let $\KK$ be a convex body in $\R^d$ and
$C_i = B_i + H_i$, $i\leq N$, be $k$-codimensional cylinders with $1\leq k<d$.
Denote $\bar C_i = C_i\cap \KK$, $i\leq N$, and $E_i=H_i^{\perp}$.
We say that the $C_i$'s form a packing in $\KK$ if $B_i\subset P_{E_i} \KK$ for
every $i\leq N$ and the interiors ${\rm int}(\bar C_i)$ of $\bar C_i$'s are pairwise disjoint.
More generally,
we say that the $C_i$'s form an $r$-fold packing in $\KK$ if $B_i\subset P_{E_i} \KK$
for every $i\leq N$ and each point of $\KK$ belongs to at most $r$ of
${\rm int}(\bar C_i)$'s. Clearly, a $1$-fold packing is just a packing.}
\end{definition}

\medskip

First we provide estimates in the case of $1$-codimensional cylinders.
Recall here that $d_{\KK}$ denotes the Banach-Mazur distance to the
Euclidean ball and $\omega _n$ denotes the volume of $\BB _2^n$.

\begin{theorem}\label{covcyl}
Let  $\KK$ be an ellipsoid in $\R^d$.
Let $C_1, \dots, C_N$ be $1$-co\-di\-men\-si\-o\-nal cylinders in $\R^d$,
which form an $r$-fold packing in $\KK$. Then
\begin{equation}
 \sum _{i=1}^N  \mbox{\rm crv} _{\KK} (C_i)
 \leq r .
\end{equation}
\end{theorem}

\medskip

\medskip

\begin{remark}
{\rm This theorem can be used to get bounds in the general case as well.
Indeed, let $\KK$ be a convex body in $\R^d$ and $T$ be an invertible linear
transformation satisfying
$$
 d_{\KK}^{-1} T\BB_2^d \subset  \KK\subset T \BB_2^d.
$$
Let $C_1, \dots, C_N$ be $1$-co\-di\-men\-si\-o\-nal cylinders in $\R^d$
forming an $r$-fold packing in $T\BB_2^d$.
Then, using definitions and Theorem~\ref{covcyl}, we observe
$$
  \sum _{i=1}^N \crv _{\KK} (C_i)
   = \sum _{i=1}^N  \frac{\vol _{d-1} (B_i)}{
  \vol _{d-1} (P_{E_i} \KK)}
  \leq \sum _{i=1}^N \frac{\vol _{d-1} (B_i)}{
  \vol _{d-1} (P_{E_i} d_{\KK}^{-1} T\BB_2^d)}
$$
$$
   =
  d_{\KK}^{d-1} \ \sum _{i=1}^N \crv _{T\BB_2^d} (C_i)
   \leq r{d_{\KK}^{d-1}}.
$$}
\end{remark}

\medskip

\proof
Every $C_i$ can be presented as $C_i = B_i + \ell _i$,
where $\ell _i$ is a line containing $0$ in $\R^d$
and $B_i$ is a body in $E_i: =\ell _i^{\perp}$ such that
$B_i\subset P_{E_i} \KK$.

Since
$\crv _{\KK} (C) = \crv _{T \KK} (T C)$ for every invertible affine map
$T: \R^d \to \R^d$, we may assume that $\KK = \BB_2^d$.  Then
$$
  \crv _{\KK} (C_i) = \frac{\vol _{d-1} (B_i)}{ \omega_{d-1}}.
$$

Consider the following (density) function on $\R^d$
$$
  p(x) = 1/\sqrt{1-|x|^2}
$$
for $|x|<1$ and $p(x)=0$ otherwise. The corresponding measure on $\R^d$
we denote by $\mu$, that is $d\mu (x) = p(x) dx$. Let $\ell$ be a line containing $0$
in $\R^d$ and $E=\ell ^{\perp}$.
It follows from direct calculations that for every $z\in E$ with $|z| < 1$
$$
  \int _{\ell +z} p(x) \ dx =  \pi.
$$
Thus we have
$$
  \mu (\BB_2^d) = \int _{\BB_2^d} p(x) \ dx = \int _{\BB_2^d \cap E}
  \int _{\ell +z} p(x) \ dx\ dz = \pi \ \omega_{d-1}
$$
and for every $i\leq N$
$$
  \mu (C_i) =\mu\left(C_i \cap \BB_2^d\r) =
  \int _{C_i} p(x)\ dx = \int _{B_i} \int _{\ell _i +z}
  p(x) \ dx\ dz = \pi \ \vol_{d-1} \left(B_i\r).
$$
Since each point of $\KK$ belongs to at most $r$ of
${\rm int}(\bar C_i)$'s, where $\bar C_i =  C_i \cap \BB_2^d$, $i\leq N$,
we obtain that
$$
 r \pi \, \omega_{d-1} =r \mu (\BB_2^d) \geq
  \sum _{i=1}^N  \mu\left(\bar C_i\r) =  \sum _{i=1}^N
  \pi \ \vol_{d-1} \left(B_i\r).
$$
This implies
\begin{equation}\label{ellcyl}
  \sum _{i=1}^N \crv _{\BB_2^d} (C_i) = \sum _{i=1}^N
  \frac{\vol _{d-1} (B_i)}{\omega_{d-1}} \leq r ,
\end{equation}
which completes the proof.
\kkk

\bigskip

Now recall the following idea from \cite{AKP}. Consider the density function defined on $\R^d$ as follows:
$p(x) = 1$ for $|x|=1$ and $p(x)=0$ otherwise. The corresponding
measure on $\R^d$ we denote by $\mu$, that is $d\mu (x) = p(x) d\lambda(x)$, where $\lambda$
is the Lebesgue measure on $S^{d-1}$. Let $H$ be a plane containing $0$ in $\R^d$ and $E= H^{\perp}$.
Then for every $z\in E$ with $|z| < 1$
$$
  \int _{H +z} p(x) \ dx =  2 \pi.
$$
(Hint: Let $\ell$ be a line parallel to $H$ and passing through $0$ in $\R^d$. Moreover, let $\ell^{\perp}= \BB_2^{d-1}$.
If $y\in(H+z)\cap S^{d-1}$ and $\alpha$ denotes the angle between the line passing through $0$ and $x=P_{\BB_2^{d-1}}y$
and the hyperplane tangent to $S^{d-1}$ at $y$, then $\cos \alpha=\sqrt{1-|x|^2}$ and so, the density at $x$ in $\BB_2^{d-1}$ is equal
to $\frac{1}{\sqrt{1-|x|^2}}$.) Therefore repeating the proof of Theorem~\ref{covcyl} with respect to the
just introduced density function $p(x)$, we obtain the following theorem.

\begin{theorem}\label{covcyl-2}
Let  $\KK$ be an ellipsoid in $\R^d$.
Let $C_1, \dots, C_N$ be $2$-co\-di\-men\-si\-o\-nal cylinders in $\R^d$,
which form an $r$-fold packing in $\KK$. Then
\begin{equation}
 \sum _{i=1}^N  \mbox{\rm crv} _{\KK} (C_i)
 \leq r.
\end{equation}
\end{theorem}

\medskip

\begin{remark}
{\rm As in the case of 1-codimensional cylinders, this theorem
can be generalized in the following way. Let $\KK$ be a convex body
in $\R^d$ and $T$ be an invertible linear transformation $T$
satisfying $d_{\KK}^{-1} T\BB_2^d \subset \KK\subset T \BB_2^d$. Let
$C_1, \dots, C_N$ be $2$-co\-di\-men\-si\-o\-nal cylinders in $\R^d$
forming an $r$-fold packing in $T\BB_2^d$. Then
\begin{equation}
 \sum _{i=1}^N  \mbox{\rm crv} _{\KK} (C_i)
 \leq r{d_{\KK}^{d-2}}.
\end{equation}}
\end{remark}

\medskip

The next theorem deals with $k$-codimensional convex cylinders.

\begin{theorem}\label{covcylgen}
Let $0<k<d$ and  $\KK$ be a convex body in $\R^d$.
Let $C_i=B_i + H_i$, $i\leq N$, be $k$-codimensional cylinders in $\R^d$,
which form an $r$-fold packing in $\KK$ and let $\bar C_i = C_i \cap \KK$.
Assume that $\bar C_i$'s  are convex bodies in $\R^d$.
Then
$$
    \sum _{i=1}^N  \mbox{\rm crv} _{\KK} (C_i) \leq
  r {d \choose k}\, \max _{i\leq N} \, \frac{\max _{x\in \R^d}\vol _k
   (\KK\cap (x+H_i)) }{\max _{x\in \R^d}\vol _k (\bar C_i\cap (x+ H_i)}.
$$
\end{theorem}

\medskip

\proof
As before denote $E_i=H_i^{\perp}$ and $B_i= P_{E_i} C_i$.
As $\bar C_i$'s form a packing in $\KK$ we have
$B_i \subset P_{E_i} \KK$ and hence $P_{E_i} \bar C_i = B_i$.

Applying Theorem~\ref{rsh} and remark following it,
we obtain for every $1\leq i\leq N$
$$
  \mbox{\rm crv} _{\KK} (C_i) =
  \frac{\vol _{d-k} (B_i)}{\vol _{d-k} (P_{E_i} \KK)}  =
  \frac{\vol _{d-k} (P_{E_i} \bar C_i)}{\vol _{d-k} (P_{E_i} \KK)}
$$
$$
  \leq {d\choose k}\, \frac{\vol _{d} (\bar C_i) }{\max _{x\in \R^d}\
  \vol _{k}  \left( \bar C_i \cap \left( x + H _i \r) \r)  } \
  \frac{\max _{x\in \R^d}\  \vol _{k}  \left(\KK \cap
  \left( x + H_i \r) \r) }{\vol _{d} (\KK)} .
$$
Since the $C_i$'s form an $r$-fold packing in $\KK$, we observe that
$$
  \sum _{i=1}^N \vol _{d} \left( \bar C_i \r) \leq r \, \vol _{d}
  \left(\KK \r) ,
$$
which implies the desired result.
\kkk


\bigskip

Finally we show an example showing some restrictions on the upper bound. We need the following simple lemma.

\begin{lemma}\label{integral}
For every $\delta \in (0, \pi/2)$ and every $n\geq 1$ one has
$$
  \frac{\delta  \, (\sin \delta)^n}{e (n+1)} \leq
  \int _{\pi/2-\delta}^{\pi/2} (\cos t)^n \, dt
  \leq \delta  \, (\sin \delta)^n.
$$
\end{lemma}

\medskip

\proof
The upper estimate is trivial, as $\cos (\cdot) $ is a decreasing function on $(0,\pi/2)$.
For the lower bound note that $\sin (\beta \delta) \geq \beta \sin \delta$ for
$\beta \in (0,1)$ and therefore
$$
  \int _{\pi/2-\delta}^{\pi/2}  (\cos t)^n \, dt \geq
   \int _{\pi/2-\delta}^{\pi/2-\beta \delta}  (\cos t)^n \, dt \geq
  (1-\beta) \delta \, (\sin (\beta \delta))^{n} \geq
  (1-\beta) \delta \beta ^{n} (\sin  \delta)^{n}.
$$
The choice $\beta = n/(n+1)$ completes the proof.
\kkk

The next theorem shows that sum of $\mbox{\rm crv} _{\KK} (C_i)$ cannot
be too small in general. The example is based on a packing of cylinders,
whose bases are caps in the Euclidean ball. We will use the following notation.
Given $m$-dimensional subspace $E$
in $\R^d$, $\delta \in (0, \pi/2)$ and $x\in S^{d-1}\cap E$ we denote
$$
  S(x, \delta, E) := \{z \in \BB_2^d \cap E \, \, \mid\, \, |\la z, x \ra |
  \geq \cos \delta\}.
$$
In other words, $S(x, \delta, E)$ is a (solid) cap in the Euclidean ball in  $E$
with the center at $x$ and the (geodesic) radius $\delta$. We also
denote
$$
   S(x, \delta ) := \{z \in S^{d-1} \, \, \mid\, \, |\la z, x \ra |
  \geq \cos \delta\},
$$
that is $S(x, \delta )$ is a (spherical) cap in the Euclidean ball in  $\R^d$.

\begin{theorem}\label{example}
Let $d>3$, $1\leq k<d$ and $\delta \in (0, \pi/4)$.
There exist  $k$-co\-di\-men\-si\-o\-nal cylinders $C_1, \dots, C_N$
in $\R^d$, which form a packing in $\BB_2^d$ and satisfy
$$
    \sum _{i=1}^N  \mbox{\rm crv} _{\BB_2^d} (C_i)=\frac{1}{\omega_{d-k}}\sum _{i=1}^N \vol _{d-k} (B_i)
   \geq
   \frac{c\, \sqrt{d}\,  (\sin \delta)^{2-k}}{2^{d-2}\, (d-k)^{3/2}} ,
$$
where  $c$ is an absolute positive constant.
\end{theorem}

\medskip

\begin{remark}
{\rm The proof of Theorem~\ref{example} below uses
representations of caps $S(x, \delta)$ as $k$-codimensional
cylinders $C(x)= S(x, \delta, E_x) + E_x^{\perp}$, where $E_x$
is a $(d-k)$-dimensional subspace $E_x$ containing $x$. Then
$\bar C(x) = C(x) \cap \BB_2^d = \conv S(x, \delta)$ and it is not
difficult to see that the maximum over $y\in \R^d$ of
$\vol _{k} (\bar C(x) \cap(y+E_x^{\perp}))$ attains at $y=\cos (\delta)x$
and  equals
$$
   \vol _{k} (y+\sqrt{1-|y|^2}\, \BB_2^d \cap E_x^{\perp})) =
    (\sin\delta)^k \omega _k.
$$
Thus, for such cylinders the upper bound from Theorem~\ref{covcylgen}
becomes
$$
  {d \choose k} \, (\sin\delta)^{-k}.
$$
Therefore, in this example, the ratio between the upper and lower bounds
is of the order $C(d, k) (\sin\delta)^{-2}$.}
\end{remark}

\medskip

\proof
Given $x\in S^{d-1}$ we construct a $k$-codimensional cylinder $C(x)$
in the following way. Fix a $(d-k)$-dimensional subspace $E_x$ containing
$x$. Let $C(x)= S(x, \delta, E_x) + E_x^{\perp}$. Of course $C(x)$ depends
on the choice of $E_x$, so for every $x$ we fix one such $E_x$. With such a
construction we have
$$
  \bar C(x) = C(x) \cap \BB_2^d = S(x, \delta, \R^d ) = \conv S(x, \delta) .
$$
Note that using the Fubini theorem and substitution $x=\sin t$, one has
$$
  \vol _{d-k} (S(x, \delta, E)) = \int _{\cos \delta}^1 (1-x^2)^{\frac{d-k-1}{2}}
  \omega _{d-k-1}\, dx
  =\omega _{d-k-1} \int_{\pi/2 - \delta}^{\pi/2} (\cos t)^{d-k} \, dt
$$
and similarly
$
  \omega _d
  =\omega _{d-1} \int_{-\pi/2}^{\pi/2} (\cos t)^{d} \, dt .
$

Now we construct a packing of caps in $S^{d-1}$ and estimate its cardinality
using standard volumetric argument. Let $\{x_i\}_{i\leq N}$ be a maximal
(in the sense of  inclusion) $(2\delta)$-separated set, that is the geodesic
distance between $x_i$ and $x_j$ is larger than $2\delta$ whenever $i\ne j$.
Then clearly the caps $S(x_i, \delta)$ are pairwise disjoint, hence so are
$S(x_i, \delta, \R^d )$'s. Therefore, the cylinders $C(x_i)$ form a packing
of $\BB_2^d$. On the other hand, due to the maximality of the set
$\{x_i\}_{i\leq N}$, the caps $S(x_i, 2\delta)$ cover $S^{d-1}$.
Therefore
$$
  1\leq \sum _{i=1}^N \sigma (S(x_i, 2\delta)) = N \sigma(S(x_1, 2\delta)),
$$
where $\sigma$ is the normalized Lebesgue measure of the sphere $S^{d-1}$.
Thus $N\geq ( \sigma(S(x_1, 2\delta)))^{-1}$. The measure of a spherical
cap can be directly calculated as
(see e.g. Chapter 2 of \cite{MS})
$$
   \sigma(S(x_1, 2\delta)) = \frac{\int_{\pi/2 - 2\delta}^{\pi/2} (\cos t)^{d-2} \, dt}
   {\int_{-\pi/2 }^{\pi/2} (\cos t)^{d-2}\, dt  } =
   \frac{\omega _{d-3}}{\omega _{d-2}}\, \int_{\pi/2 - 2\delta}^{\pi/2} (\cos t)^{d-2} \, dt .
$$
%
Finally we obtain that there are $N$ cylinders $C(x_i)$, which form packing and
$$
  A:= \sum _{i=1}^N  \mbox{\rm crv} _{\BB_2^d} (C_i) = N \,
   \frac{\vol _{d-k} (S(x, \delta, E))}{\omega _{d-k}}  \geq
   \frac{\omega _{d-k-1}}{\omega_{d-k}}\,  \frac{\omega _{d-2}}{\omega_{d-3}}\,
   \frac{  \int_{\pi/2 - \delta}^{\pi/2} (\cos t)^{d-k} \, dt }
  { \int_{\pi/2 - 2\delta}^{\pi/2} (\cos t)^{d-2} \, dt }.
$$
Using Lemma~\ref{integral} and estimates for the volume of the Euclidean ball, we
observe for an absolute positive constant $c$,
$$
  A\geq \frac{c\, \sqrt{d}}{(d-k)^{3/2} } \, \frac{(\sin \delta)^{d-k}}
  {(\sin(2\delta))^{d-2} }
  \geq   \frac{c\, \sqrt{d}\,  (\sin \delta)^{2-k}}{2^{d-2}\, (d-k)^{3/2}}.
$$
\kkk

\section{More on packing by
cylinders}
\label{sumvolumes}

In this section we estimate the total volume of bases of $1$-codimensional cylinders forming a multiple packing in a given convex body.







\begin{theorem}\label{pack1cyl}
Let  $\KK$ be a convex body in $\R^d$. For $i\leq N$ let
$C_i=B_i +H_i$ be $1$-codimensional cylinders in $\R^d$,
which form an $r$-fold packing in $\KK$. Then
\begin{equation}\label{pack-vol}
  \sum _{i=1}^N  \vol _{d-1}(B_i) \leq c_d \, r \,
  \max _{\dim L= d-1} \vol _{d-1} (P_L \KK),
\end{equation}
where
$c_d = d\, \omega _d/(2 \omega _{d-1}) \sim \sqrt{\pi \, d/ 2}$
(as $d$ grows to infinity).
\end{theorem}

\proof
We use Cauchy formula for the surface area of $\KK$:
$$
  s(\KK) = \frac{1}{\omega _{d-1}} \int _{S^{d-1}} \vol _{d-1}
  \left(P _{u^{\perp}} \KK\r) d\lambda (u) ,
$$
where $d\lambda (\cdot)$ is the Lebesgue measure on $S^{d-1}$.


For $i\leq N$ denote $\bar C_i=C_i\cap \KK$. As $C_i$'s form an $r$-fold packing in $\KK$
we have that $(\bar C_i\cap {\rm bd}\KK)$'s form an $r$-fold packing on the boundary ${\rm bd}\KK$ of $\KK$ and therefore
$$
  \sum _{i=1}^N s(\bar C_i\cap {\rm bd}\KK) \leq r\cdot s(K).
$$


Using that $\vol _{d-1}( B_i)\leq \frac{1}{2} \, s(\bar C_i\cap {\rm bd}\KK)$ and that
$\lambda (S^{d-1}) = d \omega _d$, we obtain

$$
2\sum _{i=1}^N  \vol _{d-1}(B_i) \leq r\cdot s(\KK) = \frac{r}{\omega _{d-1}} \int _{S^{d-1}} \vol _{d-1}
  \left(P _{u^{\perp}} \KK\r) d\lambda (u)
  $$
  $$
  \leq r\frac{d\, \omega _d}{ \omega _{d-1}} \,
  \max _{\dim L= d-1} \vol _{d-1} (P_L \KK),
  $$
  finishing the proof of Theorem~\ref{pack1cyl}.
\kkk

\bigskip

\begin{remark}
{\rm It would be interesting to find the best possible value of $c_d$ (as a function of $d$)
in Theorem~\ref{pack1cyl}.
Note that Theorem~\ref{covcyl} implies that when $\KK$ is an ellipsoid one can take
$c_d=1$. This leads to another natural problem:
Provide a characterization of convex bodies in $\R^d$ that satisfy
Theorem~\ref{pack1cyl} with $c_d$ bounded by an absolute constant.}
\end{remark}

\section{Packing counterpart of Falconer's bounds}
\label{Falconer}

In this section we provide a packing counterpart of Falconer's bounds.
Recall that Falconer (\cite{Fa}) gave an elegant analytic proof of the following multiple covering version of Bang's theorem
in $\R^2$ and  $\R^3$. Let $\KK$ be a convex body in $\R^2$ or $\R^3$
(i.e., a convex domain)
whose minimal width ${\rm w}(\KK)$ is equal
to the diameter of its incircle (note that any $0$-symmetric convex domain has this property). If finitely many planks form
an $r$-fold covering of $\KK$, then the sum of the widths of the planks is at least $r\, {\rm w}(\KK)$. Here we provide the packing counterpart
of Falconer's estimate in $\R^2$. Following Hadwiger (\cite{Ha}) we say, that a finite family of closed circular disks form a separable arrangement
in $\R^2$ if there exists a line in $\R^2$ that is disjoint from all the disks and divides the plane into two open half-planes each containing
at least one disk. In the opposite case we shall call the family a non-separable arrangement (in short, an NS-family) in $\R^2$. In other words,
a finite family of closed circular disks form a non-separable arrangement (i.e., an NS-family) in $\R^2$
if no line of $\R^2$ divides the disks into two non-empty sets without touching or intersecting at least one disk. We call the convex hull (resp., the sum
of the diameters) of an NS-family of disks an NS-domain (NS-diameter). The following theorem improves Theorem~\ref{pack1cyl} (with
$c_2=1$) for NS-domains $\KK$ whose NS-diameter ${\rm diam}_{NS}(\KK)$ satisfies ${\rm diam}_{NS}(\KK)=2R_{\KK}={\rm diam}({\KK})$, where $R_{\KK}$
denotes the radius of the smallest circular disk containing $\KK$ (also called the circumradius of $\KK$) and ${\rm diam}({\KK})$ denotes the Euclidean diameter of $\KK$.

\begin{theorem}\label{dual-Falconer}
Let $\KK$ be an arbitrary NS-domain in $\R^2$. If finitely many planks form an $r$-fold packing in
$\KK$, then the sum of the widths of the planks is at most $r{\rm diam}_{NS}(\KK)$. Here, $2R_{\KK}\le {\rm diam}_{NS}(\KK)$ with equality
if and only if $2R_{\KK}= {\rm diam}_{NS}(\KK)={\rm diam}({\KK})$.
\end{theorem}

Our proof is a packing analogue of Falconer's analytic method introduced for coverings by planks in \cite{Fa}. The core part of the discussions
that follow is in $\R^d$ and might be of independent interest.
Let $\KK$ be a convex body in  $\R ^d$ and let
$$\mathcal{L}^+(\KK)=\{f: \KK\to \R^+\ |\ f\ge 0\ {\rm and\ Lebesgue\ integrable\ over}\  \KK\}.$$
Moreover, let $H(s, u)=\{ x\in\R^d \ |\ \la x , u \ra=s\}$ denote the hyperplane  in $\R ^d$ with normal vector $u\in S^{d-1}$ lying at distance $s\ge 0$ from the origin $0$.
Furthermore, let the sectional integral of $f$ over $H(s,u)\cap{\rm int}(\KK)$ be denoted by
$$
  F(f, s,u)=\int_{H(s,u)\cap{\rm int}(\KK)}f(x)\ d_{H(s,u)}x
$$
for any $H(s,u)$ with $H(s,u)\cap {\rm int}(\KK) \ne\emptyset$ and with respect to the corresponding $(d-1)$-dimensional Lebesgue measure over $H(s,u)$.  Moreover, for $\Delta >0$ let
$$
  \mathcal{L}^+_{\Delta}(\KK)=\{f\in\mathcal{L}^+(\KK)\ |\ F(f, s,u)
  \ge\Delta\ {\rm for\ all\ }H(s,u)\ {\rm with}\ H(s,u)\cap {\rm int}(\KK) \ne\emptyset\}.
$$
Finally, let
$$
  m(\mathcal{L}^+_{\Delta}(\KK))=\inf\left\{\int_{\KK}f(x)dx\ |\ f\in \mathcal{L}^+_{\Delta}(\KK)\right\}.
$$
If $g$ is a Lebesgue integrable function on $\R$
then $g(\la x, u\ra)$, $x\in \R^d$, is called a ridge
function in the direction $u\in S^{d-1}$.

\begin{lemma}\label{ridge-estimate}
Let $\KK$ be a convex body in $\R^d$ and let $u_i\in S^{d-1}$, $1\le i\le N$.
Let $g_i(\la x, u_i\ra)$, $1\le i\le N$, be ridge functions such that the support of $g_i$ is contained in $[a_i, b_i]$, where
$$
  a_i=\min\{\la x, u_i\ra\ |\ x\in\KK\} \, \, \, \mbox{ and }
  \, \, \, b_i=\max\{\la x, u_i\ra\ |\ x\in\KK\}.
$$
Assume that for every $x\in\KK$,
$$
  \sum_{i=1}^N g_i(\la x, u_i\ra)\le 1 .
$$
Then
$$
\sum_{i=1}^N \int_{-\infty}^{+\infty}g_i(t)dt\le m(\mathcal{L}^+_1(\KK)).
$$
\end{lemma}

\proof
For every $f\in\mathcal{L}^+_1(\KK)$ one has
$$
  \sum_{i=1}^N \int_{-\infty}^{+\infty}g_i(t)dt\le \sum_{i=1}^N \int_{-\infty}^{+\infty}g_i(t)F(f, t, u_i)dt
$$
$$
  =\int_{\KK}f(x)\sum_{i=1}^N g_i(\la x, u_i\ra)dx\le\int_{\KK}f(x)dx ,
$$
which implies the desired result.
\kkk

\medskip

Let $\KK$ be a convex body in $\R ^d$ and let    $\BB(\KK)=x_{\KK}+R_{\KK}\BB_2^d$ denote the
circumscribed ball of $\KK$, i.e. the smallest Euclidean ball
containing $\KK$ with center at $x_{\KK}$ and radius $R_{\KK}$,
called the circumradius of $\KK$.
Recall that the support function of $\KK$ is defined by
$$
 h_{\KK}(x)=\sup\{\la x, k\ra\ |\ k\in \KK\}
$$
for $x\in\R^d$.

\begin{lemma}\label{inf-estimate}
If $\KK$ is a convex body with circumradius $R_{\KK}$ in $\R^d$, then
$$m(\mathcal{L}^+_{\Delta}(\KK))\ge 2\Delta R_{\KK}.$$
\end{lemma}

\proof
Let $f\in \mathcal{L}^+_{\Delta}(\KK)$. As $\int_{\KK}f(x)dx>0$,
$\KK$ weighted by $f$ has a centroid, i.e., a point $c$ such that
$$
 \int_{\KK}f(x)(x-c)dx=0,
$$
which lies inside $\KK$. Passing to the body $\KK_c=\KK-c$, without loss
of generality we may assume that $c$ is the origin $0$.
Let $h_{\KK}$ be the support function of $\KK$ and $R$ be the smallest number such that $\KK\subset R\BB^d_2$. Then $R\geq R_{\KK}$ and
there exists  $u\in S^{d-1}$ with
\begin{equation}\label{direction-u}
h_{\KK}(u) = R \geq R_{\KK}.
\end{equation}
Taking moments perpendicular to $u$ yield
\begin{equation}\label{moments}
\int_{-h_{\KK}(-u)}^{h_{\KK}(u)}tF(f,t,u)dt=0.
\end{equation}
Using $F(f,t,u)\ge\Delta$ and (\ref{direction-u}) we observe
\begin{equation}\label{first-lowerbound}
\int_{0}^{h_{\KK}(u)}tF(f,t,u)dt\ge\frac{1}{2}\Delta R_{\KK}^2.
\end{equation}
Then (\ref{moments}) and (\ref{first-lowerbound}) yield that
\begin{equation}\label{second-lowerbound}
\int_{0}^{h_{\KK}(-u)}tF(f,t,u)dt\ge\frac{1}{2}\Delta R_{\KK}^2.
\end{equation}
One can check that
\begin{equation}\label{third-bound}
\inf\left\{\int_0^{A}F(t)dt\ |\ A>0, \, F(t)\ge\Delta\, \, {\rm  and }\, \,  \int_0^{+\infty}tF(t)dt\ge M\right\}=(2M\Delta)^{\frac{1}{2}}
\end{equation}
with the infimum being attained if and only if
$F(t) = \Delta$ for (almost)
all $t\le A$ and $A=(\frac{2M}{\Delta})^{\frac{1}{2}}$. Thus (\ref{first-lowerbound}), (\ref{second-lowerbound}), and (\ref{third-bound}) yield
$$\int_{\KK}f(x)dx=\int_{-h_{\KK}(-u)}^{h_{\KK}(u)}F(f,t,u)dt\ge 2\left(2\left(\frac{1}{2}\Delta R_{\KK}^2\right)\Delta\right)^{\frac{1}{2}}=2\Delta R_{\KK},$$
finishing the proof of Lemma~\ref{inf-estimate}.
\kkk

\medskip

\noindent
{\bf Proof of Theorem~\ref{dual-Falconer}:\, }
Let $\KK$ be an NS-domain in $\R^2$ with NS-diameter ${\rm diam}_{NS}(\KK)$
and let $C_1, \dots , C_N$ be planks that form an $r$-fold packing in $\KK$.
For every $1\le i\leq N$, choose $u_i\in S^{d-1}$ which is orthogonal $C_i$ and
let
$$
  a_i=\min\{\la x, u_i\ra\ |\ x\in C_i \} \, \, \, \mbox{ and }
  \, \, \, b_i=\max\{\la x, u_i\ra\ |\ x\in C_i \}.
$$
Clearly, the Euclidean width $w(C_i)$ of the plank $C_i$ satisfies
$w(C_i)= b_i-a_i$. Consider the ridge functions
$$
   g_i(\langle x, u_i \rangle)  =
    \frac{1}{r} \, {\chi}_{[a_i, b_i]}(\langle x, u_i \rangle),
$$
where ${\chi}_{[a_i, b_i]}$ is the characteristic function of the
segment $[a_i, b_i]$.
On the one hand, Lemma~\ref{ridge-estimate} applied to $g_i$'s implies that
\begin{equation}\label{I}
\sum_{i=1}^{N} w(C_i)\le rm(\mathcal{L}^+_1(\KK)).
\end{equation}
On the other hand, recall the following well-known fact (Theorem 1.2 in \cite{Fa}):
$m(\mathcal{L}^+_1(R\BB_2^d))=2R$ and this value is attained uniquely by the function $f(x)=\frac{1}{\pi R}(R^2-|x|^2)^{-\frac{1}{2}}$ if $|x|<R$ and $f(x)=0$ if $|x|\ge R$.
By taking the sum of the analogue functions over the generating circular disks of the NS-domain $\KK$ we get that
\begin{equation}\label{II}
m(\mathcal{L}^+_1(\KK))\le {\rm diam}_{NS}(\KK).
\end{equation}
This and (\ref{I}) yield $\sum_{i=1}^{N} w(C_i)\le r {\rm diam}_{NS}(\KK)$. Moreover, Lemma~\ref{inf-estimate}  and (\ref{II})  imply
that
\begin{equation}\label{III}
2R_{\KK}\le {\rm diam}_{NS}(\KK).
\end{equation}
 For a completely different proof of (\ref{III}) we refer
the interested reader to Goodman and Goodman (\cite{GoGo}).
As the case of equality is rather straightforward to show,
this completes the proof of Theorem~\ref{dual-Falconer}.
\kkk

%





\bigskip

\noindent
{\bf Acknowledgement.} The authors would like to thank R. Karasev
for comments on Remark~\ref{multiple-covering-by-planks} and
the anonymous referee for remarks and careful reading.

\vspace{1cm}

\medskip

\noindent
K\'aroly Bezdek,
Department of Mathematics and Statistics,
2500 University drive N.W.,
University of Calgary, AB, Canada, T2N 1N4.
\newline
{\sf e-mail: bezdek@math.ucalgary.ca}

\smallskip

\noindent
A.E. Litvak,
Department  of Mathematical and Statistical Sciences,
University of Alberta, Edmonton, AB, Canada, T6G 2G1.
\newline
{\sf e-mail: aelitvak@gmail.com}


\begin{thebibliography}{GGM}



\bibitem[AKP]{AKP} A.~Akopyan, R.~Karasev, and F.~Petrov, {\em Bang's problem
and symplectic invariants}, arXiv:1404.0871v1 (2014), 1--15.


\bibitem[B1]{Ball} K.~Ball,
{\em  Flavors of geometry} in
{\em An elementary introduction to modern convex geometry},
Levy, Silvio (ed.), Cambridge:
Cambridge University Press. Math. Sci. Res. Inst. Publ. 31, 1--58 (1997).




\bibitem[B2]{Ball2} K.~Ball, {\em The plank problem for symmetric bodies},
 Invent. Math. 104 (1991),  535--543.


\bibitem[Ba]{Ba}T.~Bang, {\em A solution of the ``plank problem"},
Proc. Amer. Math. Soc.  2 (1951), 990--993.











\bibitem[BL]{BL} K.~Bezdek, A.E.~Litvak, {\em Covering convex bodies by
cylinders and lattice points by flats},
J. Geom. Analysis, 19 (2009), 233--243.


\bibitem[C]{Ch} G. D. Chakerian,  {\em Inequalities for the difference body of
a convex body},  Proc. Amer. Math. Soc.  18 (1967), 879--884.

\bibitem[Fa]{Fa}K. J. Falconer, {\em Function space topologies defined by sectional integrals
and applications to an extremal problem}, Math. Proc. Camb. Phil. Soc. 87 (1980), 81--96.

\bibitem[Fe]{Fe}L. Fejes T\'oth, {\em Lagerungen in der Ebene, auf der Kugel und im Raum},
Grundlehren Math. Wiss. 65, Springer, Berlin, 1953.

\bibitem[GG]{GoGo}A. W. Goodman and R. E. Goodman, {\em A circle covering theorem},
Amer. Math. Monthly 52 (1945), 494--498.

\bibitem[H]{Ha}H. Hadwiger, {\em Nonseparable convex systems}, Amer. Math. Monthly 54 (1947), 583--585.








\bibitem[J]{J} F.~John,
{\em Extremum problems with inequalities as subsidiary conditions},
Studies and Essays Presented to R. Courant on his 60th Birthday,
January 8, 1948, 187--204. Interscience Publishers, Inc., New York,
N. Y., 1948.












\bibitem[MS]{MS} V. D. Milman, G. Schechtman,
{\em Asymptotic theory of finite-dimensional
normed spaces.} With an appendix by M. Gromov,
Lect. Notes in Math., 1200.
Springer-Verlag, Berlin, 1986.



\bibitem[Pi]{Pi2} G.~Pisier,
{\em The Volume of Convex Bodies and Banach Space Geometry},
Cambridge University Press 1989.


\bibitem[RS]{RS} C. A. Rogers,  G. C. Shephard,
{\em Convex bodies associated with a given convex body},
J. London Math. Soc.  33 (1958), 270--281.








\end{thebibliography}
\end{document}